# Non-dominated Solution of Fuzzy Maximum-Return Problem

U. M. Pirzada

*School of Science and Engineering , Navrachana University*
*Vaodara-391410, India*
Email: salmapirzada@yahoo.com

D. C. Vakaskar

*Department of Applied Mathematics, Faculty of Tech. and Engg.,*
*M. S. University of Baroda,, Vadodara-390001, Gujarat, India.*
Email: dcvakaskar@gmail.com

Abstract:

In this paper, we find a non-dominated solution of a fuzzy maximum-return problem (unconstrained single-variable fuzzy optimization problem). We establish Newton method to find the solution of the unconstrained single-variable fuzzy optimization problem using the differentiability of $\alpha$-level functions of a fuzzy-valued function and partial order relation on a set of fuzzy numbers.

Key words:
Fuzzy-valued functions, Newton method, Maximum-return problem

**2000 MSC:** 03E72, 90C70

## 1. Introduction

Maximum-return problem is formulated in [2]. He has proposed an optimal portfolio by considering a composite function of return and risk using weight parameter. He has minimized the composite function using some single-variable non-linear optimization methods. The value of acceptable risk is known approximately in the problem. So, rather than to consider the approximate value of acceptable risk as a exact number, we consider it as a fuzzy number, which capture the real meaning of approximation. By this way, we fuzzify the crisp optimization to study the flexibility in the approximate data.

---

Received April, 2016

☆ The work is supported by National Board for Higher Mathematics (NBHM), Department of Atomic Energy (DAE), India.

* Corresponding author. Tel.: +91 9825979539

Fuzzy optimization accounts for any imprecision in the optimization problem. Bellman and Zadeh (1970) introduced fuzzy optimization problems in [3] where they have stated that a fuzzy decision can be viewed as the intersection of fuzzy goals and problem constraints. Afterwards, a lot of articles dealing with linear and non-linear fuzzy optimization problems were published. The optimality conditions for fuzzy optimization problems have been proposed by several authors ( See [10], [11] and [12]).

Newton method is one of the most used method in numerical analysis to find the root of equations. The method is well-established in crisp optimization theory for solving unconstrained and constrained optimization problems and it is area of interest for many researchers (ref. [1, 15]). It has numerous applications in management science, industrial and financial research, data mining. The advantage of the Newton method is its speed of convergence once a sufficiently accurate approximation is known. To the best of our knowledge, iterative methods are not explored for fuzzy optimization problems. In [14], authors have proposed the Newton method for solving multi-variable fuzzy optimization problem. Under the concept of Hukuhara differentiability ($H$-differentiability) of a fuzzy-valued function due to Puri and Ralescu [16], the method is proposed. The $H$-differentiability is the most used differentiability to study the fuzzy differential systems. Unfortunately, only limited class of fuzzy-valued functions are $H$-differentiable, since this differentiability is based on existence of Hukuhara difference. For instance, (See [5])

Let $\tilde{f}:(0,2\pi) \to F(\mathbb{R})$ be defined on $\alpha$-level sets by

$$\left[\tilde{f}(x)\right]_\alpha = (1-\alpha)(2+\sin(x))[-1,1], \tag{1}$$

for $\alpha \in [0,1]$. The function is not $H$-differentiable at $x^0 = \pi/2$.

For these kinds of functions, we can not apply the Newton method proposed in [14]. This motivates us to develop the Newton method under the differentiability of $\alpha$-level functions of a fuzzy-valued function.

In this paper, we consider an unconstrained single-variable fuzzy optimization problem with respect to a specific partial order relation defined on set of fuzzy numbers. Under the concept of the differentiability of $\alpha$-level functions of a fuzzy-valued function over $R$, we propose the Newton method to find the non-dominated solution of the problem. One illustration is given to explain the application of the proposed method. We also consider the fuzzy maximum-return problem as a case study which studies the effect of the parameter at its exact value as well as the numbers nears to it with varying degree of membership.

The paper is organized in following manner. In section 2, we cite some basic definitions and results on fuzzy numbers, partial order relation on fuzzy numbers and differentiability of $\alpha$-level functions of a fuzzy-valued function over $R$. We present the Newton method to find a non-dominated solution of an unconstrained single-variable fuzzy optimization problem in Section 3. Convergence of the proposed method is studied in the same section. Non-dominated solution of fuzzy maximum-return problem is given as a case study in Section 4. Moreover, one more illustration is given as a application of the proposed method. We conclude in the last section.

**2. Preliminaries**

We start with some basic definitions.

**Definition 2.1** [6]. Let $\mathbb{R}$ be the set of real numbers and $\tilde{a}:\mathbb{R}\to[0,1]$ be a fuzzy set. We say that $\tilde{a}$ is a fuzzy number if it satisfies the following properties:

(i) $\tilde{a}$ is normal, that is, there exists $x_0 \in \mathbb{R}$ such that $\tilde{a}(x_0)=1$;

(ii) $\tilde{a}$ is fuzzy convex, that is, $\tilde{a}(tx+(1-t)y)\geq \min\{\tilde{a}(x),\tilde{a}(y)\}$, whenever $x,y \in \mathbb{R}$ and $t \in [0,1]$;

(iii) $\tilde{a}(x)$ is upper semi-continuous on $\mathbb{R}$, that is, $\{x|\tilde{a}(x)\geq \alpha\}$ is a closed subset of $\mathbb{R}$ for each $\alpha \in (0,1]$;

(iv) $cl\{x \in \mathbb{R}|\tilde{a}(x)>0\}$ forms a compact set,

where $cl$ denotes closure of a set. The set of all fuzzy numbers on $\mathbb{R}$ is denoted by $F(\mathbb{R})$. For all $\alpha \in (0,1]$, $\alpha$-level set $\tilde{a}_\alpha$ of any $\tilde{a} \in F(\mathbb{R})$ is defined as $\tilde{a}_\alpha = \{x \in \mathbb{R}|\tilde{a}(x)\geq \alpha\}$. The 0-level set $\tilde{a}_0$ is defined as the closure of the set $\{x \in \mathbb{R}|\tilde{a}(x)>0\}$. By definition of fuzzy numbers, we can prove that, for any $\tilde{a} \in F(\mathbb{R})$ and for each $\alpha \in (0,1]$, $\tilde{a}_\alpha$ is compact convex subset of $\mathbb{R}$, and we write $\tilde{a}_\alpha = [\tilde{a}_\alpha^L, \tilde{a}_\alpha^U]$. $\tilde{a} \in F(\mathbb{R})$ can be recovered from its $\alpha$-level sets by a well known decomposition theorem (ref. [7]), which states that $\tilde{a} = \bigcup_{\alpha \in [0,1]} \alpha \cdot \tilde{a}_\alpha$ where union on the right-hand side is the standard fuzzy union.

**Definition 2.2.** According to Zadeh's extension principle, we define addition and scalar multiplication as well as multiplication and division of two fuzzy numbers $\tilde{a}, \tilde{b} \in F(\mathbb{R})$ by their $\alpha$-level sets as follows:

$$(\tilde{a}\oplus \tilde{b})_\alpha = [\tilde{a}_\alpha^L + \tilde{b}_\alpha^L, \tilde{a}_\alpha^U + \tilde{b}_\alpha^U]$$

$$(\lambda \odot \tilde{a})_\alpha = [\lambda \cdot \tilde{a}_\alpha^L, \lambda \cdot \tilde{a}_\alpha^U], \text{ if } \lambda \geq 0$$

$$= [\lambda \cdot \tilde{a}_\alpha^U, \lambda \cdot \tilde{a}_\alpha^L], \text{ if } \lambda < 0, \ \lambda \in \mathbb{R}$$

$$(\tilde{a}\otimes \tilde{b})_\alpha = \left[\min\{\tilde{a}_\alpha^L \cdot \tilde{b}_\alpha^L, \tilde{a}_\alpha^L \cdot \tilde{b}_\alpha^U, \tilde{a}_\alpha^U \cdot \tilde{b}_\alpha^L, \tilde{a}_\alpha^U \cdot \tilde{b}_\alpha^U\}, \max\{\tilde{a}_\alpha^L \cdot \tilde{b}_\alpha^L, \tilde{a}_\alpha^L \cdot \tilde{b}_\alpha^U, \tilde{a}_\alpha^U \cdot \tilde{b}_\alpha^L, \tilde{a}_\alpha^U \cdot \tilde{b}_\alpha^U\}\right]$$

$$(\tilde{a}\oslash \tilde{b})_\alpha = \left[\min\{\tilde{a}_\alpha^L / \tilde{b}_\alpha^L, \tilde{a}_\alpha^L / \tilde{b}_\alpha^U, \tilde{a}_\alpha^U / \tilde{b}_\alpha^L, \tilde{a}_\alpha^U / \tilde{b}_\alpha^U\}, \max\{\tilde{a}_\alpha^L / \tilde{b}_\alpha^L, \tilde{a}_\alpha^L / \tilde{b}_\alpha^U, \tilde{a}_\alpha^U / \tilde{b}_\alpha^L, \tilde{a}_\alpha^U / \tilde{b}_\alpha^U\}\right]$$

where level sets of $\tilde{a}$ and $\tilde{b}$ are $\tilde{a}_\alpha = [\tilde{a}_\alpha^L, \tilde{a}_\alpha^U]$, $\tilde{b}_\alpha = [\tilde{b}_\alpha^L, \tilde{b}_\alpha^U]$ for $\alpha \in [0,1]$. Moreover, $0 \notin [\tilde{b}_\alpha^L, \tilde{b}_\alpha^U]$ for division operation.

**Definition 2.3** [19]. The membership function of a triangular fuzzy number $\tilde{a}$ is defined as

$$\tilde{a}(r) = \begin{cases} (r-a^L)/(a-a^L) & \text{if } a^L \leq r \leq a \\ (a^U - r)/(a^U - a) & \text{if } a < r \leq a^U \\ 0 & \text{otherwise} \end{cases}$$

which is denoted by $\tilde{a} = (a^L, a, a^U)$. Then $\alpha$-level set of $\tilde{a}$ is then

$$\tilde{a}_\alpha = \left[(1-\alpha)a^L + \alpha a, (1-\alpha)a^U + \alpha a\right].$$

Based on multiplication and division operations on fuzzy numbers, we define the square and multiplicative inverse of a triangular fuzzy number as follows:

**Definition 2.4.** Let $\tilde{a} = (a^L, a, a^U)$ be a triangular fuzzy number. We assume that $a^L, a, a^U \geq 0$. The square and multiplicative inverse of $\tilde{a}$ using its $\alpha$-level sets are defined as

$$(\tilde{a}^2)_\alpha = \left[\left((1-\alpha)a^L + \alpha a\right)^2, \left((1-\alpha)a^U + \alpha a\right)^2\right]$$

$$(1/\tilde{a})_\alpha = \left[1/\left((1-\alpha)a^U + \alpha a\right), 1/\left((1-\alpha)a^L + \alpha a\right)\right],$$

where $0 \notin \left[(1-\alpha)a^L + \alpha a, (1-\alpha)a^U + \alpha a\right]$, for all $\alpha \in [0,1]$ to define multiplicative inverse operation.

Among the various ordering methods for fuzzy numbers the commonly used one is a partial order relation called fuzzy-max order, introduced by Ramik and Rimanek [17], which is defined as follows:

**Definition 2.5.** Let $\tilde{a}$ and $\tilde{b}$ be two fuzzy numbers in $F(\mathbb{R})$ and let $\tilde{a}_\alpha = \left[\tilde{a}_\alpha^L, \tilde{a}_\alpha^U\right]$ and $\tilde{b}_\alpha = \left[\tilde{b}_\alpha^L, \tilde{b}_\alpha^U\right]$ be two closed intervals in $\mathbb{R}$, $\alpha \in [0,1]$. We define $\tilde{a} \preceq \tilde{b}$ if and only if $\tilde{a}_\alpha^L \leq \tilde{b}_\alpha^L$ and $\tilde{a}_\alpha^U \leq \tilde{b}_\alpha^U$ for all $\alpha \in [0,1]$. Moreover, $\tilde{a} \prec \tilde{b}$ if and only if $\tilde{a} \preceq \tilde{b}$ and there exists an $\alpha_0 \in [0,1]$ such that $\tilde{a}_{\alpha_0}^L < \tilde{b}_{\alpha_0}^L$ or $\tilde{a}_{\alpha_0}^U \leq \tilde{b}_{\alpha_0}^U$.

We define comparable fuzzy numbers as follows:

**Definition 2.6.** For $\tilde{a}$, $\tilde{b}$ in $F(\mathbb{R})$, we say that $\tilde{a}$ and $\tilde{b}$ are comparable if either $\tilde{a} \preceq \tilde{b}$ or $\tilde{b} \preceq \tilde{a}$.

**Definition 2.7** [8]. Let $V$ be a real vector space and $F(\mathbb{R})$ be a set of fuzzy numbers. Then a function $\tilde{f}: V \to F(\mathbb{R})$ is called fuzzy-valued function defined on $V$. Corresponding to such a function $\tilde{f}$ and $\alpha \in [0,1]$, we define two real-valued functions $\tilde{f}_\alpha^L$ and $\tilde{f}_\alpha^U$ on $V$ as $\tilde{f}_\alpha^L(x) = \left(\tilde{f}(x)\right)_\alpha^L$ and $\tilde{f}_\alpha^U(x) = \left(\tilde{f}(x)\right)_\alpha^U$ for all $x \in V$. These functions $\tilde{f}_\alpha^L(x)$ and $\tilde{f}_\alpha^U(x)$ are called $\alpha$-level functions of the fuzzy-valued function $\tilde{f}$.

The Hausdorff metric on a set of fuzzy numbers is defined as follows:

**Definition 2.8** [9]. Let $A, B \subseteq \mathbb{R}^n$. The Hausdorff metric $d_H$ is defined by

$$d_H(A,B) = \max\left\{\sup_{x \in A}\inf_{y \in B}\|x-y\|, \sup_{y \in B}\inf_{x \in A}\|x-y\|\right\}.$$

Then the metric $d_F$ on $F(\mathbb{R})$ is defined as

$$d_F(\tilde{a},\tilde{b}) = \sup_{0 \leq \alpha \leq 1}\left\{d_H(\tilde{a}_\alpha,\tilde{b}_\alpha)\right\},$$

for all $\tilde{a}, \tilde{b} \in F(\mathbb{R})$. Since $\tilde{a}_\alpha$ and $\tilde{b}_\alpha$ are closed bounded intervals in $\mathbb{R}$,

$$d_F(\tilde{a},\tilde{b}) = \sup_{0 \leq \alpha \leq 1}\max\left\{\left|\tilde{a}_\alpha^L - \tilde{b}_\alpha^L\right|, \left|\tilde{a}_\alpha^U - \tilde{b}_\alpha^U\right|\right\}.$$

Pirzada and Pathak [14] have studied Newton method using Hukuhara differentiability of a fuzzy-valued function. This differentiability is based on Hukuhara difference. Hukuhara difference between two fuzzy numbers define as follows:

**Definition 2.9.** Let $\tilde{a}$ and $\tilde{b}$ be two fuzzy numbers. If there exists a fuzzy number $\tilde{c}$ such that $\tilde{c} \oplus \tilde{b} = \tilde{a}$. Then $\tilde{c}$ is called Hukuhara difference of $\tilde{a}$ and $\tilde{b}$ be is denoted by $\tilde{a} \ominus_H \tilde{b}$.

We have the following remark (ref. [5]):

**Remark 2.1.** The necessary condition for the Hukuhara difference $\tilde{a} \ominus_H \tilde{b}$ to exist is that some translate of $\tilde{b}$ is a fuzzy subset of $\tilde{a}$. For instance, $\tilde{a} = (-1,1,3)$ and $\tilde{b} = (-1,0,1)$ are triangular fuzzy numbers such that $(-1,1,3) = (0,1,2) \oplus (-1,0,1)$ then $(0,1,2)$ is called $H$-difference of $\tilde{a}$ and $\tilde{b}$, as $\tilde{b}$ is a fuzzy subset of $\tilde{b}$.

Hukuhara differentiability ($H$-differentiability) of a fuzzy-valued function due to Puri and Ralescu [16] define as follows:

**Definition 2.10.** Let $(a,b) \subset \mathbb{R}$. A fuzzy-valued function $\tilde{f}:(a,b) \to F(\mathbb{R})$ is said to be $H$-differentiable at $x_0 \in (a,b)$ if there exists a fuzzy number $\tilde{f}'(x_0)$ such that the limits (with respect to $d_F$)

$$\lim_{h \to 0^+}\frac{1}{h} \odot \left[\tilde{f}(x_0+h) \ominus_H \tilde{f}(x_0)\right], \text{ and } \lim_{h \to 0^+}\frac{1}{h} \odot \left[\tilde{f}(x_0) \ominus_H \tilde{f}(x_0-h)\right]$$

both exist and are equal to $\tilde{f}'(x_0)$. In this case, $\tilde{f}'(x_0)$ is called the $H$-derivative of $\tilde{f}$ at $x_0$. If $\tilde{f}$ is $H$-differentiable at any $x \in [a,b]$, we can $\tilde{f}$ is $H$-differentiable on $[a,b]$.

Existence of $H$-differentiability of a fuzzy-valued function depend upon the existence of Hukuhara difference. For instance, a fuzzy-valued function (1) given in previous section is not $H$-differentiable. However, under some conditions Hukuhara difference exists and therefore, we say that function $H$-differentiable. The conditions are given the following proposition.

**Proposition 2.1** [4]. *Let $\tilde{c} \in F(\mathbb{R})$ and $g:(a,b) \to \mathbb{R}_+$ be differentiable at $x_0 \in (a,b)$. Define $\tilde{f}:(a,b) \to F(\mathbb{R})$ by $\tilde{f}(x) = \tilde{c} \odot g(x)$, for all $x \in (a,b)$. If we suppose that, $g'(x_0) > 0$ then Hukuhara differences of $\tilde{f}$ exist and $\tilde{f}$ is $H$-differentiable at $x_0$ with $\tilde{f}'(x) = \tilde{c} \odot g'(x)$.*

The limited class of fuzzy-valued functions satisfies the above conditions. To overcome this situation, in this paper we use the differentiability of $\alpha$-level functions $\tilde{f}_\alpha^L(x)$ and $\tilde{f}_\alpha^U(x)$, all $\alpha \in [0,1]$ of a fuzzy-valued function $\tilde{f}$ to proposed the Newton method. The $\alpha$-level functions $\tilde{f}_\alpha^L(x)$ and $\tilde{f}_\alpha^U(x)$, all $\alpha \in [0,1]$ are real-valued functions.

## 3. Newton method

We consider an unconstrained single variable fuzzy optimization problem $(USFOP)$:

$$\text{Minimize } \tilde{f}(x), \ x \in X$$

Where $X \subseteq \mathbb{R}$ is an open set and $\tilde{f}: X \to F(\mathbb{R})$ is a fuzzy-valued function.

Using partial order relation defined in Section 2, we define a locally non-dominated solution of $(USFOP)$ as follows:

**Definition 3.1.** Let $X \subseteq \mathbb{R}$ be an open set and let $\tilde{f}: X \to F(\mathbb{R})$ be a fuzzy-valued function. A point $x^0 \in X$ is said to be a locally non-dominated solution, if there exists no $x^1 \in N_\varepsilon(x^0) \cap X$ such that $\tilde{f}(x^1) \preceq \tilde{f}(x^0)$, where $N_\varepsilon(x^0)$ is $\varepsilon$-neighbourhood of $x^0$.

To propose the Newton method for solving an unconstrained single variable fuzzy optimization, we prove the following theorem.

**Theorem 3.1.** *Let $\tilde{f}: X \to F(\mathbb{R})$ be a fuzzy-valued function, $X$ is an open subset of $\mathbb{R}$. If $\overline{x}^0 \in X$ is a locally non-dominated solution of $(USFOP)$ and for any direction* *d and for any $\delta > 0$ there exists $\lambda \in ]0,\delta[$ such that $\tilde{f}(x^0 + \lambda \cdot d)$ and $\tilde{f}(x^0)$ are comparable, then $x^0$ is a local minimizer of $\tilde{f}_\alpha^L(x^0)$ and $\tilde{f}_\alpha^U(x^0)$, for all $\alpha \in [0,1]$.*

Proof. Suppose our claim is not true. That is, $x^0$ is not a local minimizer of $\tilde{f}_\alpha^L(x^0)$ or $\tilde{f}_\alpha^U(x^0)$, for any $\alpha_0 \in [0,1]$. Without loss of any generality, we say that $\dot{x}^0$ is not a local minimizer of $\tilde{f}_{\alpha_0}^L(x^0)$. That is

$$\tilde{f}_{\alpha_0}^L(x^0) > \tilde{f}_{\alpha_0}^L(x)$$

for at least one $x$ in every $N_\varepsilon(x^0)$.

We can always choose $x$ as $x = x^0 + \lambda \cdot d$, for $\lambda \in ]0,\delta[$, where $\delta > 0$ and $d$ is any direction such that

$$\tilde{f}_{\alpha_0}^L(x^0) > \tilde{f}_{\alpha_0}^L(x^0 + \lambda \cdot d) \tag{2}$$

Now by assumption of the theorem, for any direction $d$ and for any $\delta > 0$ there exists $\lambda \in ]0,\delta[$ such that $\tilde{f}(x^0 + \lambda \cdot d)$ and $\tilde{f}(x^0)$ are comparable. Thus, either $\tilde{f}(x^0 + \lambda \cdot d) \preceq \tilde{f}(x^0)$ or $\tilde{f}(x^0) \preceq \tilde{f}(x^0 + \lambda \cdot d)$. But from (2), we must have

$$\tilde{f}(x^0 + \lambda \cdot d) \prec \tilde{f}(x^0)$$

This inequality contradicts to our assumption that $x^0$ is a non-dominiated solution. Hence the proof.

Now we propose the Newton method to find the non-dominated solution of an unconstrained single variable fuzzy optimization problem. We assume that at each measurement point $x^{(k)}$ we can calculate $\tilde{f}_\alpha^L(x^{(k)})$, $\tilde{f}_\alpha^U(x^{(k)})$, $\left(\tilde{f}_\alpha^L(\tilde{x}^{(k)})\right)'$, $\left(\tilde{f}_\alpha^U(\tilde{x}^{(k)})\right)'$, $\left(\tilde{f}_\alpha^L(\tilde{x}^{(k)})\right)''$ and $\left(\tilde{f}_\alpha^U(\tilde{x}^{(k)})\right)''$ for all $\alpha \in [0,1]$. Therefore, we can approximate $\tilde{f}_\alpha^L$ and $\tilde{f}_\alpha^U$ by quadratic real-valued functions $\tilde{q}_\alpha^L$ and $\tilde{q}_\alpha^U$ respectively, for all $\alpha \in [0,1]$. Usign Taylor's formula

$$\tilde{q}_\alpha^L(x) = \tilde{f}_\alpha^L(x^{(k)}) + \left\{\left(\tilde{f}_\alpha^L(x^{(k)})\right)' \cdot (x - x^{(k)})\right\} + \left\{\left(\tilde{f}_\alpha^L(x^{(k)})\right)'' \cdot \left((x - x^{(k)})^2 / 2!\right)\right\}$$

and

$$\tilde{q}_\alpha^U(x) = \tilde{f}_\alpha^U(x^{(k)}) + \left\{ \left(\tilde{f}_\alpha^U(x^{(k)})\right)' \cdot (x - x^{(k)}) \right\} + \left\{ \left(\tilde{f}_\alpha^U(x^{(k)})\right)'' \left((x - x^{(k)})^2 / 2!\right) \right\}$$

for all $\alpha \in [0,1]$. If $x$ is a non-dominated solution of $(USFOP)$ then by Theorem 3.1,

$$\tilde{f}'(x) = \tilde{0}.$$

This implies

$$\left(\tilde{f}_\alpha^L(x)\right)' = \left(\tilde{f}_\alpha^U(x)\right)' = 0,$$

for all $\alpha \in [0,1]$. We get the first order necessary conditions for $\tilde{q}_\alpha^L$ and $\tilde{q}_\alpha^U$:

$$\left(\tilde{q}_\alpha^L(x)\right)' = \left(\tilde{q}_\alpha^U(x)\right)' = 0,$$

for all $\alpha \in [0,1]$. That is

$$\left(\tilde{f}_\alpha^L(x^{(k)})\right)' + \left(\tilde{f}_\alpha^L(x^{(k)})\right)'' \cdot (x - x^{(k)}) = 0$$

and

$$\left(\tilde{f}_\alpha^U(x^{(k)})\right)' + \left(\tilde{f}_\alpha^U(x^{(k)})\right)'' \cdot (x - x^{(k)}) = 0,$$

for all $\alpha \in [0,1]$. That is,

$$\int_0^1 \left(\tilde{f}_\alpha^L(x^{(k)})\right)' d\alpha + \int_0^1 \left(\tilde{f}_\alpha^L(x^{(k)})\right)'' d\alpha \cdot (x - x^{(k)}) = 0 \qquad (2)$$

and

$$\int_0^1 \left(\tilde{f}_\alpha^U(x^{(k)})\right)' d\alpha + \int_0^1 \left(\tilde{f}_\alpha^U(x^{(k)})\right)'' d\alpha \cdot (x - x^{(k)}) = 0 \qquad (3)$$

By adding (2) and (3), we have

$$\int_0^1 \left\{ \left(\tilde{f}_\alpha^L(x^{(k)})\right)' + \left(\tilde{f}_\alpha^U(x^{(k)})\right)' \right\} d\alpha + \int_0^1 \left\{ \left(\tilde{f}_\alpha^L(x^{(k)})\right)'' + \left(\tilde{f}_\alpha^U(x^{(k)})\right)'' \right\} d\alpha \cdot (x - x^{(k)}) = 0.$$

We define a real-valued function $F$ in following way.

$$F(x) = \int_0^1 \left\{ \tilde{f}_\alpha^L(x) + \tilde{f}_\alpha^L(x) \right\} d\alpha$$

Therefore, we can write the above equation as

$$F'(x^{(k)}) + F''(x^{(k)})(x - x^{(k)}) = 0 \qquad (4)$$

where

$$F'\left(x^{(k)}\right) = \int_0^1 \left\{ \left(\tilde{f}_\alpha^L\left(x^{(k)}\right)\right)' + \left(\tilde{f}_\alpha^U\left(x^{(k)}\right)\right)' \right\} d\alpha$$

and

$$F''\left(x^{(k)}\right) = \int_0^1 \left\{ \left(\tilde{f}_\alpha^L\left(x^{(k)}\right)\right)'' + \left(\tilde{f}_\alpha^U\left(x^{(k)}\right)\right)'' \right\} d\alpha .$$

By putting $x = x^{(k+1)}$ in (4), we get

$$x^{(k+1)} = x^{(k)} - \left(F'\left(x^k\right) / F''\left(x^k\right)\right) \tag{5}$$

Thus starting with an initial approximation to minimizer of $\tilde{f}$, we can generate a sequence of approximations to the minimizer of $\tilde{f}$ using the formula (5). The procedure is terminated when $\left|x^{(k+1)} - x^{(k)}\right| < \varepsilon$, $\varepsilon$ is pre-specified positive real number.

**Remark 3.1.** The method is well-defined only when $F''\left(x^{(k)}\right) \neq 0$ for each $k$.

To prove the convergence of Newton method, we need following Theorem from Numerical analysis.

**Theorem 3.2** [18]. *Let $x = x^*$ be a root of $f(x) = 0$ and let $I$ be an interval containing the point $x = x^*$. Let $\phi(x)$ and $\phi'(x)$ be continuous on $I$, where $\phi(x)$ is defined by the equation $x = \phi(x)$ which is equivalent to $f(x) = 0$. Then if $\left|\phi'(x)\right| < 1$ for all $x \in I$, the sequence of approximations $x_0, x_1, \cdots, x_n$ defined by*

$$x_{n+1} = \phi(x_n)$$

*converges to the root $x^*$, provided that the initial approximation $x_0$ is chosen in $I$.*

Now we show the convergence of Newton method.

**Theorem 3.3.** *Suppose that $\tilde{f}_\alpha^L(x)$ and $\tilde{f}_\alpha^U(x)$, for all $\alpha \in [0,1]$ are three times continuously differentiable functions defined on $\mathbb{R}$ and $x^* \in \mathbb{R}$ a point such that*
(1) $F'\left(x^*\right) = 0$
(2) $F''\left(x^*\right) \neq 0$.

*Then for all $x^0$ sufficiently close to $x^*$, Newton method is well-defined for all $x$ and converges to $x^*$ with order of convergence at least 2.*

*Proof.* Since $\tilde{f}_\alpha^L(x)$ and $\tilde{f}_\alpha^U(x)$, for all $\alpha \in [0,1]$ are three times continuously differentiable, $F''(x)$ is continuous function. We have,

$$|F''(x)| \geq \varepsilon, \tag{6}$$

for some $\varepsilon > 0$ in a suitable neighbourhood of $x^*$. Within this neighbourhood we can select an interval $I$ such that, for all $x \in I$

$$|F'(x)F'''(x)| < \varepsilon^2, \tag{7}$$

this is possible since $F'(x^*) = 0$ and $F(x)$ is also three times continuously differentiable function.

Now taking

$$\phi(x) = x - \left(F'(x)/F''(x)\right)$$

We observe that

$$\phi'(x) = F'(x)F'''(x)/\left(F''(x)\right)^2,$$

for all $x \in I$. Hence, from (6) and (7), we get. $|\phi'(x)| < 1$.

Thus, $\phi(x)$ satisfying all the hypothesis of Theorem 3.2, by taking an initial approximation $x^{(0)} \in I$, we get the sequence of approximations $x^{(0)}, x^{(1)}, \cdots, x^{(n)}$ satisfying (5) converge to non-dominated solution $x = x^*$.

Now to obtain the rate of convergence of Newton method, we note that $F'(x^*) = 0$ so that Taylor's expansion gives

$$F'\left(x^{(k)}\right) + F''\left(x^{(k)}\right)\left(x^* - x^{(k)}\right) + F'''\left(x^{(k)}\right)\left(\left(x^* - x^{(k)}\right)^2/2!\right) + \cdots = 0,$$

where

$$F'\left(x^{(k)}\right) = \int_0^1 \left\{\left(\tilde{f}_\alpha^L\left(x^{(k)}\right)\right)' + \left(\tilde{f}_\alpha^U\left(x^{(k)}\right)\right)'\right\} d\alpha,$$

$$F''\left(x^{(k)}\right) = \int_0^1 \left\{\left(\tilde{f}_\alpha^L\left(x^{(k)}\right)\right)'' + \left(\tilde{f}_\alpha^U\left(x^{(k)}\right)\right)''\right\} d\alpha$$

and

$$F'''\left(x^{(k)}\right) = \int_0^1 \left\{\left(\tilde{f}_\alpha^L\left(x^{(k)}\right)\right)''' + \left(\tilde{f}_\alpha^U\left(x^{(k)}\right)\right)'''\right\} d\alpha.$$

From which we obtain

$$-\left(F'\left(x^{(k)}\right)/F''\left(x^{(k)}\right)\right) = \left(x^* - x^{(k)}\right) + \frac{1}{2}\left(x^* - x^{(k)}\right)^2 \left(F'''\left(x^{(k)}\right)/F''\left(x^{(k)}\right)\right) \tag{8}$$

From (5) and (8), we have

$$x^{(k+1)} - x^* = (1/2)\left(x^{(k)} - x^*\right)^2 \left(F'''\left(x^{(k)}\right)/F''\left(x^{(k)}\right)\right) \tag{9}$$

Setting
$$\varepsilon_k = x^{(k)} - x^*,$$

Equation (9) gives
$$\varepsilon_{k+1} \propto \left(-F'''(x^*)/2F''(x^*)\right) \cdot \left(\varepsilon_{(k)}\right)^2,$$

so that the Newton's method has quadratic convergence.

## 4. Numerical application and case study-Fuzzy maximum return problem

First we present the algorithm of proposed Newton method.

**Algorithm 1** Newton method
1: Input $x^0$, $\varepsilon$
2: Workout $F'(x)$ and $F''(x)$
3: $k \leftarrow 0$
4: **repeat**
5: $x^{(k+1)} = x^{(k)} - \left(F'(x^{(k)})/F''(x^{(k)})\right)$
6: $k \leftarrow k+1$
7: **until**
8: $\left|x^{(k+1)} - x^{(k)}\right| < \varepsilon$
9: Optimal solution $x^* \leftarrow x^{(k)}$

Here we give one example to illustrate the application of the proposed method.

**Example 4.1.** Maximze $\tilde{f}(x) = \left(\tilde{1} \odot x^3\right) \oplus \left(\tilde{2} \odot x^2\right)$, $x \in \mathbb{R}_+$

Where $\tilde{1} = (0,1,2)$ and $\tilde{2} = (1,2,3)$ are triangular fuzzy numbers and initial approximation for maximizer is $x^0 = 1$.

Clearly, the level functions $\tilde{f}_\alpha^L(x) = (\alpha)x^3 + (1+\alpha)x^2$ and $\tilde{f}_\alpha^U(x) = (2-\alpha)x^3 + (3-\alpha)x^2$ are differentiable real-valued function for $x \in \mathbb{R}_+$ at each $\alpha \in [0,1]$.

Moreover, for any $x, y \in \mathbb{R}_+$ such that $x > y$, we have

$$\tilde{f}_\alpha^L(x) - \tilde{f}_\alpha^L(y) = \alpha(x^3 - y^3) + (1+\alpha)(x^2 - y^2) > 0,$$

$$\tilde{f}_\alpha^U(x) - \tilde{f}_\alpha^U(y) = (2-\alpha)(x^3 - y^3) + (3-\alpha)(x^2 - y^2) > 0,$$

for all $\alpha \in [0,1]$. That is, $\tilde{f}(y) \preceq \tilde{f}(x)$. Therefore, $\tilde{f}(x)$ and $\tilde{f}(y)$ are comparable). So we can use Newton method to solve this problem.

**Table 1:** Convergence of Newton method: Example 4.1

| $k$ | $x^{(k)}$ | $x^{(k+1)}$ | $\tilde{f}\left(x^{(k)}\right)$ |
|---|---|---|---|
| 0 | 1 | −0.6667 | $(-14,-11,-8)$ |
| 1 | −0.6667 | −0.2037 | $(-5.9753,-5.1358,-4.2963)$ |
| 2 | −0.2037 | −0.0201 | $(-0.5412,-0.4962,-0.4513)$ |
| 3 | −0.0201 | $-2.0150e-004$ | $(-0.0053,-0.0049,-0.0044)$ |
| 4 | $-2.0150e-004$ | $-2.0301e-008$ | $1.0e-006*(-0.5278,-0.4872,-0.4466)$ |
| 5 | $-2.0301e-008$ | $-2.0301e-008$ | $1.0e-014*(-0.5358,-0.4945,-0.4533)$ |

**Example 4.2.** In [2], author has proposed an optimal portfolio by considering a composite function of return and risk using weight parameter. The composite function is $F=-R+\rho V$. First term is the negative of the expected return and second term is a multiple of the risk. The maximum-return problem involves minimizing the function:

$$F(x)=-0.06667x-1.1167+\left(\rho/V_a^2\right)\left(0.1256x^2-0.1589x+0.05139-V_a\right)^2$$

where $V_a$ denotes an acceptable value for risk. The positive constant $\rho$ is weight parameter, controls the balance between return and risk. The problem is to find invested fraction $x$ which minimize $F$ to obtain a large value for return coupled with a small value for risk.

The approximate value of $V_a$ is 0.00168.

**Step 1.** Fuzzification of a crisp problem

In order to choose more flexible value of $V_a$, here we fuzzify the $V_a$. That is, we consider the approximate $V_a$ (near to $V_a$) as a triangular fuzzy number which capture the meaning of approximation.

Therefore, we define approximate $V_a$ as $\tilde{V}_a=\left(V_a^L,V_a,V_a^U\right)$. We choose the appropriate values of $V_a^L$, $V_a$, $V_a^U$, according to the physical situation.

Here we consider the problem by taking fuzzy data $\bar{V}_\alpha$ as $\bar{V}_a=(0.00167,0.00168, 0.00172)$. By considering the fuzzy approximate value, we can study not only the effect of the exact parameter value $V_a=0.000168$ but the simultaneous effect of the real numbers which are near to the $V_a$ with varying membership values in the same problem. Moreover, we take the weight parameter $\rho$ is also fuzzy number.

The fuzzy maximum-return problem is

$$\tilde{F}(x)=-0.06667x-1.1167+\left(\tilde{\rho}/\tilde{V}_a^2\right)\left(0.125x^2-0.1589x+0.05139-\tilde{V}_a\right)^2,$$

$\tilde{\rho}=(0.5,1.5,3.5)$ is a triangular fuzzy number.

**Step 2.** Non-dominated solution of the problem

By fuzzy arithmetic defined in section 2, we have the $\alpha$-level functions

**Table 2.** Solution of crisp maximum-return problem Example 4.2

| $V_a$ | $\rho$ | $x^*$ | $F(x^*)$ |
|---|---|---|---|
| 0.00168 | 1 | 0.6989 | -1.1633 |
| 0.00168 | 1.5 | 0.6988 | -1.1633 |
| 0.00169 | 1.5 | 0.6994 | -1.1633 |
| 0.00169 | 2 | 0.6993 | -1.1633 |

**Table 3.** Non-dominated solution of fuzzy maximum-return problem Example 4.2

| $\tilde{V}_a$ | $\tilde{\rho}$ | $x^*$ | $dF(x^*)$ |
|---|---|---|---|
| (0.00167,0.00168,0.00172) | (0.5,1.5,3.5) | 0.6988 | -1.1631 |

$$\tilde{F}_\alpha^L(x) = -0.06667x - 1.1167 + \left(((1-\alpha)0.5 + 1.5\alpha)/((1-\alpha)0.0017 + 0.00168\alpha)^2\right)$$
$$\left(0.1256x^2 - 0.1589x + 0.05139 - (1-\alpha)0.0017 - 0.00168\alpha\right),$$

and

$$\tilde{F}_\alpha^U(x) = -0.06667x - 1.1167 + \left(((1-\alpha)3.5 + 1.5\alpha)/((1-\alpha)0.00167 + 0.00168\alpha)^2\right)$$
$$\left(0.1256x^2 - 0.1589x + 0.05139 - (1-\alpha)0.00167 - 0.00168\alpha\right),$$

$\alpha \in [0,1]$. The non-dominated solution of the fuzzy optimization problem is obtained using the proposed Newton method by taking an initial root $x_0 = 1$.

**Step 3.** Comparison of solutions
Solution of crisp maximum-return problem problem is given Table 2. $V_a$ is acceptable risk, $\rho$ is weight parameter, $x^*$ is final root and $F(x^*)$ is value of minimum. In every case, method is converged with the accuracy of $10^{-5}$.

Non-dominated solution of fuzzy maximum-return problem problem is shown in Table 3. $dF(x^*)$ is the defuzzified value of fuzzy-valued minimum. We use the defuzzification method proposed in [20] and implemented its programming using $MATLAB$. The method is converged with the accuracy of $10^{-5}$.

## 5. Conclusions

Fuzzification of maximum-return problem is successfully studied. Solution of the problem is done using the Newton method. The method is proposed to find the non-dominated solution of fuzzy unconstrained single-variable optimization problem. The problem studied is a representative one but the method developed can be utilized for range of different applications where objective function may involve fuzzy parameters.

Moreover, we have used the differentiability of $\alpha$-level functions $\tilde{f}_\alpha^L$ and $\tilde{f}_\alpha^U$ corresponding to fuzzy-valued function $\tilde{f}$ to proposed the method. We observed that we can cover wider class of fuzzy-valued functions under this differentiability than the Hukuhara differentiability to optimize the fuzzy-valued functions under the Newton method proposed in [14].

**References**

[1] V. Antony Vijesh and P. V. Subrahmanyam, A Newton-like method and its application, *Journal of Mathematical Analysis and Applications*, 339 (2008), 1231-1242.

[2] M. Bartholomew-Biggs, Non-linear Optimization with Financial Applications. Springer, (2005).

[3] R. E. Bellman and L. A. Zadeh, Decision making in a fuzzy environment, Management Science, 17 (1970), 141-164.

[4] B. Bede and S. G. Gal, Generalizations of the differentiability of fuzzy-number-valued functions with applications to fuzzy differential equations, *Fuzzy Sets and Systems,* 151 (2005), 581-599.

[5] P. Diamond, Kloeden, Metric spaces of fuzzy sets, Theory and Applications, World Scientific, (1994).

[6] A. A. George, Fuzzy Ostrowski Type Inequalities, Computational and Applied mathematics, Vol. 22 (2003), 279-292.

[7] A. A. George, Fuzzy Taylor Formulae, *CUBO*, *A Mathematical Journal*, Vol. 7 (2005), 1-13.

[8] Wu Hsien-Chung, Duality Theory in Fuzzy Optimization Problems, Fuzzy Optimization and Decision Making. 3 (2004), 345-365.

[9] Wu Hsien-Chung, An $(\alpha,\beta)$-Optimal Solution Concept in Fuzzy Optimization Problems, *Optimization*, Vol. 53 (2004), 203-221.

[10] Wu Hsien-Chung, The Karush-Kuhn-Tucker optimality conditions for optimization Problems with fuzzy-valued objective functions, *Math Meth Oper Res*, 66 (2007), 203-224.

[11] Wu Hsien-Chung, The Optimality Conditions for Optimization Problems with fuzzy-valued objective functions, *Optimization*, 57 (2008), 473-489.

[12] V. D. Pathak and U. M. Pirzada, The optimality conditions for fuzzy optimization problem under the concept of generalized convexity, *Advances in Applied Mathematical Analysis,* 1 (2010), 23-38.

[13] U. M. Pirzada and V. D. Pathak, Newton Method for Solving Multi-variable Fuzzy Optimization Problem, *Journal of Optimization Theory and Applications*; Springer, 156 (2013), 867-881, DOI: 10.1007/s10957-012-0141-3.

[14] B. T. Polyak, Newtons method and its use in optimization, *European Journal of Operation Research*, 181 (2007), 1086-1096.

[15] M. L. Puri and D. A. Ralescu, Differentials of fuzzy functions, *J. of Math. Analysis and App,* 91 (1983), 552-558.

[16] J. Ramik and J. Rimanek, Inequality relation between fuzzy numbers and its use in fuzzy optimization, *Fuzzy Sets and Systems,* 16 (1985), 123-138.

[17] S. S. Sastry, Introductory Methods of Numerical Analysis (4th Ed.), PHI Learning Private Limited, (2010)

[18] S. Saito and H. Ishii, $L$-Fuzzy Optimization Problems by Parametric Representation, IEEE, (2001), 1173-1177.

[19] Y. -M. Wang, Centroid defuzzification and the maximizing set and minimizing set ranking based on alpha level sets, Computers & Industrial Engineering, (2009), doi: 10.1016/j.cie.2008.11.014.